# On the 0-Rotatability of Regular Spider Trees of Leg Length 4

**Abstract**

A tree $T$ on $q$ edges is said to be *0-rotatable* if for every vertex $v \in V(T)$, there exists a graceful labeling $f : V(T) \to \{0, 1, \dots, q\}$ such that $f(v) = 0$. Rofa proved that all regular (symmetric) spider trees with leg length at most 3 are 0-rotatable and conjectured that the property holds for every leg length. In this note, we settle the first open case of this conjecture by proving that every regular spider tree $T(k, 4)$ with $k \geq 1$ legs of length 4 is 0-rotatable. Our proof combines Rofa's algebraic graceful labeling at the root with complementary labelings, an end-to-end leg-reversal construction placing 0 at a leaf, and a two-leg reflection permutation placing 0 at distance 2 from the root.

## 1 Introduction

Let $G = (V, E)$ be a simple graph with $q = |E|$ edges. A *graceful labeling* of $G$, introduced by Rosa [5] and popularized by Golomb [2], is an injection $f : V(G) \to \{0, 1, \dots, q\}$ such that the induced edge weight function $w_f : E(G) \to \{1, 2, \dots, q\}$, defined by

$$w_f(uv) = |f(u) - f(v)| \quad \text{for all } uv \in E(G),$$

is a bijection. The longstanding Ringel–Kotzig conjecture asserts that every tree admits a graceful labeling.

To construct graceful labelings inductively by amalgamating smaller trees at designated root vertices, Chung and Hwang [1] introduced the notion of $k$-rotatability: a tree $T$ is *0-rotatable* if for every vertex $v \in V(T)$, there exists a graceful labeling $f$ of $T$ satisfying $f(v) = 0$. While several families of trees—including paths, caterpillars, and stars—are known to be 0-rotatable, determining 0-rotatability for general symmetric trees remains challenging because a separate graceful construction is required for each automorphism orbit of vertices.

A *regular spider tree* $T(k, \ell)$ consists of a central root vertex $v_0$ of degree $k$ to which $k$ pairwise internally disjoint paths (legs) of uniform length $\ell$ are attached. Using an algebraic graceful labeling framework for rooted symmetric trees [3], Rofa [4] established that $T(k, \ell)$ is 0-rotatable for all $k \geq 1$ whenever $\ell \leq 3$, and conjectured that $T(k, \ell)$ is 0-rotatable for every $\ell \geq 1$. In this paper, we resolve the case $\ell = 4$ completely.

## 2 Preliminaries and Notation

For any integer $k \geq 1$, let $T(k, 4)$ denote the regular spider tree with central root $v_0$ and $k$ legs $L_1, L_2, \dots, L_k$, each of length 4. For each $i \in \{1, 2, \dots, k\}$, denote the vertices along leg $L_i$ in order of increasing distance from $v_0$ by $v_{i,1}, v_{i,2}, v_{i,3}, v_{i,4}$, so that

$$V(T(k,4)) = \{v_0\} \cup \bigcup_{i=1}^{k} \{v_{i,1}, v_{i,2}, v_{i,3}, v_{i,4}\},$$

$$E(T(k,4)) = \bigcup_{i=1}^{k} \{v_0 v_{i,1}, v_{i,1} v_{i,2}, v_{i,2} v_{i,3}, v_{i,3} v_{i,4}\}.$$

Thus $T(k, 4)$ has $n = |V(T(k, 4))| = 4k + 1$ vertices and $q = |E(T(k, 4))| = 4k$ edges.

**Lemma 1** (Complementary Labeling). *If $f : V(G) \to \{0, 1, \ldots, q\}$ is a graceful labeling of a graph $G$ with $q$ edges, then the complementary labeling $\bar{f} : V(G) \to \{0, 1, \ldots, q\}$ defined by $\bar{f}(v) = q - f(v)$ is also graceful.*

*Proof.* Since $v \mapsto q - v$ is a bijection on $\{0, 1, \ldots, q\}$, $\bar{f}$ is a bijection onto $\{0, 1, \ldots, q\}$. Moreover, for every edge $uv \in E(G)$,

$$w_{\bar{f}}(uv) = |\bar{f}(u) - \bar{f}(v)| = |(q - f(u)) - (q - f(v))| = |f(u) - f(v)| = w_f(uv),$$

so $w_{\bar{f}}$ induces the identical bijection onto $\{1, 2, \ldots, q\}$. □

Because the automorphism group $\mathrm{Aut}(T(k, 4))$ acts transitively on the set of legs $\{L_1, \ldots, L_k\}$, the vertex set $V(T(k, 4))$ partitions into exactly five automorphism orbits indexed by distance from the root $v_0$:

$$O_0 = \{v_0\} \quad \text{and} \quad O_d = \{v_{1,d}, v_{2,d}, \ldots, v_{k,d}\} \quad \text{for } d \in \{1, 2, 3, 4\}.$$

Consequently, $T(k, 4)$ is 0-rotatable if and only if there exist five graceful labelings $f_0, f_1, f_2, f_3, f_4$ satisfying $f_0(v_0) = 0$ and $f_d(v_{1,d}) = 0$ for each $d \in \{1, 2, 3, 4\}$.

## 3 Main Result

**Theorem 2.** *For every integer $k \geq 1$, the regular spider tree $T(k, 4)$ is 0-rotatable.*

*Proof.* Let $q = 4k$. We construct a graceful labeling placing 0 in each of the five vertex orbits $O_0, O_1, O_4, O_3, O_2$ in turn.

***Case 0: 0 at the Root ($O_0$).*** Since $T(k, 4)$ is a rooted symmetric tree of height 4, Rofa's algebraic graceful labeling function [3, 4] directly provides a graceful labeling $f_0 : V(T(k, 4)) \to \{0, 1, \ldots, 4k\}$ defined by $f_0(v_0) = 0$ and, for each $i \in \{1, 2, \ldots, k\}$,

$$f_0(v_{i,1}) = 4(k - i + 1), \qquad f_0(v_{i,2}) = 4i - 3, \qquad f_0(v_{i,3}) = 4(k - i + 1) - 1, \qquad f_0(v_{i,4}) = 4i - 2.$$

In particular, $f_0$ is graceful with $f_0(v_0) = 0$ and $f_0(v_{1,1}) = 4k$.

***Case 1: 0 at Distance 1 ($O_1$).*** Define $f_1(v) = 4k - f_0(v)$ to be the complementary labeling of $f_0$. By Lemma 1, $f_1$ is graceful, and

$$f_1(v_{1,1}) = 4k - f_0(v_{1,1}) = 4k - 4k = 0.$$

***Case 4: 0 at a Leaf ($O_4$).*** Define $f_4 : V(T(k, 4)) \to \{0, 1, \ldots, 4k\}$ piecewise by setting $f_4(v_0) = 2$,

$$f_4(v_{1,1}) = 4k - 1, \qquad f_4(v_{1,2}) = 1, \qquad f_4(v_{1,3}) = 4k, \qquad f_4(v_{1,4}) = 0,$$

and, for each $i \in \{2, 3, \ldots, k\}$,

$$f_4(v_{i,1}) = 4(k - i + 1) + 2, \qquad f_4(v_{i,2}) = 4i - 5, \qquad f_4(v_{i,3}) = 4(k - i + 1) + 1, \qquad f_4(v_{i,4}) = 4i - 4.$$

We verify that $f_4$ is a graceful labeling:

(i) *Vertex bijectivity.* Partitioning the image of $f_4$ by residue classes modulo 4 gives

$$\{f_4(v_{1,4})\} \cup \{f_4(v_{i,4})\}_{i=2}^{k} \cup \{f_4(v_{1,3})\} = \{0\} \cup \{4, 8, \ldots, 4k-4\} \cup \{4k\},$$
$$\{f_4(v_{1,2})\} \cup \{f_4(v_{i,3})\}_{i=2}^{k} = \{1\} \cup \{5, 9, \ldots, 4k-3\},$$
$$\{f_4(v_0)\} \cup \{f_4(v_{i,1})\}_{i=2}^{k} = \{2\} \cup \{6, 10, \ldots, 4k-2\},$$
$$\{f_4(v_{i,2})\}_{i=2}^{k} \cup \{f_4(v_{1,1})\} = \{3, 7, \ldots, 4k-5\} \cup \{4k-1\},$$

which are pairwise disjoint and whose union is $\{0, 1, \ldots, 4k\}$.

(ii) *Edge bijectivity.* Along leg $L_1$, the four edges induce weights

$$w_{f_4}(v_0v_{1,1}) = 4k - 3, \qquad w_{f_4}(v_{1,1}v_{1,2}) = 4k - 2, \qquad w_{f_4}(v_{1,2}v_{1,3}) = 4k - 1, \qquad w_{f_4}(v_{1,3}v_{1,4}) = 4k.$$

For each $i \in \{2, \ldots, k\}$, the root-incident edge receives weight $w_{f_4}(v_0v_{i,1}) = 4(k - i + 1)$, covering $\{4, 8, \ldots, 4k - 4\}$. Meanwhile, the three internal edges along $L_i$ ($2 \le i \le k$) induce the consecutive weight triple

$$w_{f_4}(v_{i,1}v_{i,2}) = |4k - 8i + 11|, \qquad w_{f_4}(v_{i,2}v_{i,3}) = |4k - 8i + 10|, \qquad w_{f_4}(v_{i,3}v_{i,4}) = |4k - 8i + 9|.$$

Re-indexing via $j = i - 1 \in \{1, \ldots, k - 1\}$, these triples become $(|4k - 8j + 3|, |4k - 8j + 2|, |4k - 8j + 1|)$. As $j$ increases from 1 to $\lceil (k - 1)/2 \rceil$, the middle term $4k - 8j + 2$ decreases in steps of 8 from $4k - 6$ down to $4k - 8\lceil (k - 1)/2 \rceil + 2 \in \{2, 6\}$; as $j$ increases from $\lceil (k - 1)/2 \rceil + 1$ to $k - 1$, the middle term $|4k - 8j + 2| = 8j - 4k - 2$ increases in steps of 8 from $8\lceil (k - 1)/2 \rceil - 4k + 6 \in \{6, 2\}$ up to $4k - 10$. Thus these triples cover every integer in $\{1, 2, \ldots, 4k - 5\}$ not divisible by 4. Combining all edges yields $w_{f_4}(E) = \{1, 2, \ldots, 4k\}$, so $f_4$ is graceful with $f_4(v_{1,4}) = 0$ and $f_4(v_{1,3}) = 4k$.

***Case 3: 0 at Distance 3 ($O_3$).*** Define $f_3(v) = 4k - f_4(v)$ to be the complementary labeling of $f_4$. By Lemma 1, $f_3$ is graceful, and

$$f_3(v_{1,3}) = 4k - f_4(v_{1,3}) = 4k - 4k = 0.$$

***Case 2: 0 at Distance 2 ($O_2$).*** We consider two subcases according to $k$:

- *Subcase* $k = 1$. Here $T(1, 4)$ is the path $v_0v_{1,1}v_{1,2}v_{1,3}v_{1,4}$ with $q = 4$. Setting

  $$(f_2(v_0), f_2(v_{1,1}), f_2(v_{1,2}), f_2(v_{1,3}), f_2(v_{1,4})) = (2, 3, 0, 4, 1)$$

  produces edge weights 1, 3, 4, 2 in order along the path, so $f_2$ is graceful with $f_2(v_{1,2}) = 0$.

- *Subcase* $k \ge 2$. Define $f_2 : V(T(k, 4)) \to \{0, 1, \ldots, 4k\}$ by setting $f_2(v) = f_4(v)$ for all $v \in \{v_0\} \cup \bigcup_{i=3}^{k} V(L_i)$ and modifying the labels on $L_1 \cup L_2$ as follows:

  $$f_2(v_{1,1}) = 4k - 2, \qquad f_2(v_{1,2}) = 0, \qquad f_2(v_{1,3}) = 4k, \qquad f_2(v_{1,4}) = 1,$$
  $$f_2(v_{2,1}) = 4k - 1, \qquad f_2(v_{2,2}) = 4, \qquad f_2(v_{2,3}) = 4k - 3, \qquad f_2(v_{2,4}) = 3,$$

  and, for each $i \in \{3, 4, \ldots, k\}$ (when $k \ge 3$),

  $$f_2(v_{i,1}) = 4(k - i + 1) + 2, \qquad f_2(v_{i,2}) = 4i - 5, \qquad f_2(v_{i,3}) = 4(k - i + 1) + 1, \qquad f_2(v_{i,4}) = 4i - 4.$$

  Since $k \ge 2$, we have $4k - 3 \ge 5 > 4$, so the eight values assigned to $V(L_1) \cup V(L_2)$ under $f_2$ form the exact same set $\{0, 1, 3, 4, 4k - 3, 4k - 2, 4k - 1, 4k\}$ as under $f_4$. Therefore $f_2$ is a bijection onto $\{0, 1, \ldots, 4k\}$.

  The eight edges incident to $L_1 \cup L_2$ realize the weights

  $$w_{f_2}(v_0v_{1,1}) = 4k - 4, \qquad w_{f_2}(v_{1,1}v_{1,2}) = 4k - 2, \qquad w_{f_2}(v_{1,2}v_{1,3}) = 4k, \qquad w_{f_2}(v_{1,3}v_{1,4}) = 4k - 1,$$
  $$w_{f_2}(v_0v_{2,1}) = 4k - 3, \qquad w_{f_2}(v_{2,1}v_{2,2}) = 4k - 5, \qquad w_{f_2}(v_{2,2}v_{2,3}) = 4k - 7,$$
  $$w_{f_2}(v_{2,3}v_{2,4}) = 4k - 6,$$

whose union is $\{4k, 4k - 1, \ldots, 4k - 7\}$. Under $f_4$, the eight edges incident to $L_1 \cup L_2$ realized $\{4k, 4k - 1, 4k - 2, 4k - 3\} \cup \{4k - 4, 4k - 5, 4k - 6, 4k - 7\}$, the identical set of weights. Since $f_2$ agrees with $f_4$ on $\{v_0\} \cup \bigcup_{i=3}^{k} V(L_i)$, the remaining $4k - 8$ edges realize $\{1, 2, \ldots, 4k - 8\}$. Hence $f_2$ is graceful with $f_2(v_{1,2}) = 0$.

In all five cases, every orbit $O_d$ $(0 \le d \le 4)$ admits a graceful labeling assigning 0 to a representative vertex. Thus $T(k, 4)$ is 0-rotatable for all $k \ge 1$. □

## 4 Illustration for $T(5, 4)$

To illustrate Theorem 2, Table 1 and Figure 1 display the explicit vertex sequences and induced edge weights for the regular spider $T(5, 4)$ ($k = 5$, $q = 20$) across the three primary constructions $f_0$, $f_4$, and $f_2$ (with $f_1 = 20 - f_0$ and $f_3 = 20 - f_4$ obtained immediately by complementation).

Table 1: Graceful labelings $f_0$, $f_4$, and $f_2$ of $T(5, 4)$ ($q = 20$).

| **Labeling** | **Leg** | **Vertex Labels ($v_0 - v_{i,1} - v_{i,2} - v_{i,3} - v_{i,4}$)** | **Induced Edge Weights** |
|---|---|---|---|
| $f_0$ ($0 \in O_0$) | $L_1$ | **0** − 20 − 1 − 19 − 2 | 20, 19, 18, 17 |
| | $L_2$ | **0** − 16 − 5 − 15 − 6 | 16, 11, 10, 9 |
| | $L_3$ | **0** − 12 − 9 − 11 − 10 | 12, 3, 2, 1 |
| | $L_4$ | **0** − 8 − 13 − 7 − 14 | 8, 5, 6, 7 |
| | $L_5$ | **0** − 4 − 17 − 3 − 18 | 4, 13, 14, 15 |
| $f_4$ ($0 \in O_4$) | $L_1$ | 2 − 19 − 1 − 20 − **0** | 17, 18, 19, 20 |
| | $L_2$ | 2 − 18 − 3 − 17 − 4 | 16, 15, 14, 13 |
| | $L_3$ | 2 − 14 − 7 − 13 − 8 | 12, 7, 6, 5 |
| | $L_4$ | 2 − 10 − 11 − 9 − 12 | 8, 1, 2, 3 |
| | $L_5$ | 2 − 6 − 15 − 5 − 16 | 4, 9, 10, 11 |
| $f_2$ ($0 \in O_2$) | $L_1$ | 2 − 18 − **0** − 20 − 1 | 16, 18, 20, 19 |
| | $L_2$ | 2 − 19 − 4 − 17 − 3 | 17, 15, 13, 14 |
| | $L_3$ | 2 − 14 − 7 − 13 − 8 | 12, 7, 6, 5 |
| | $L_4$ | 2 − 10 − 11 − 9 − 12 | 8, 1, 2, 3 |
| | $L_5$ | 2 − 6 − 15 − 5 − 16 | 4, 9, 10, 11 |

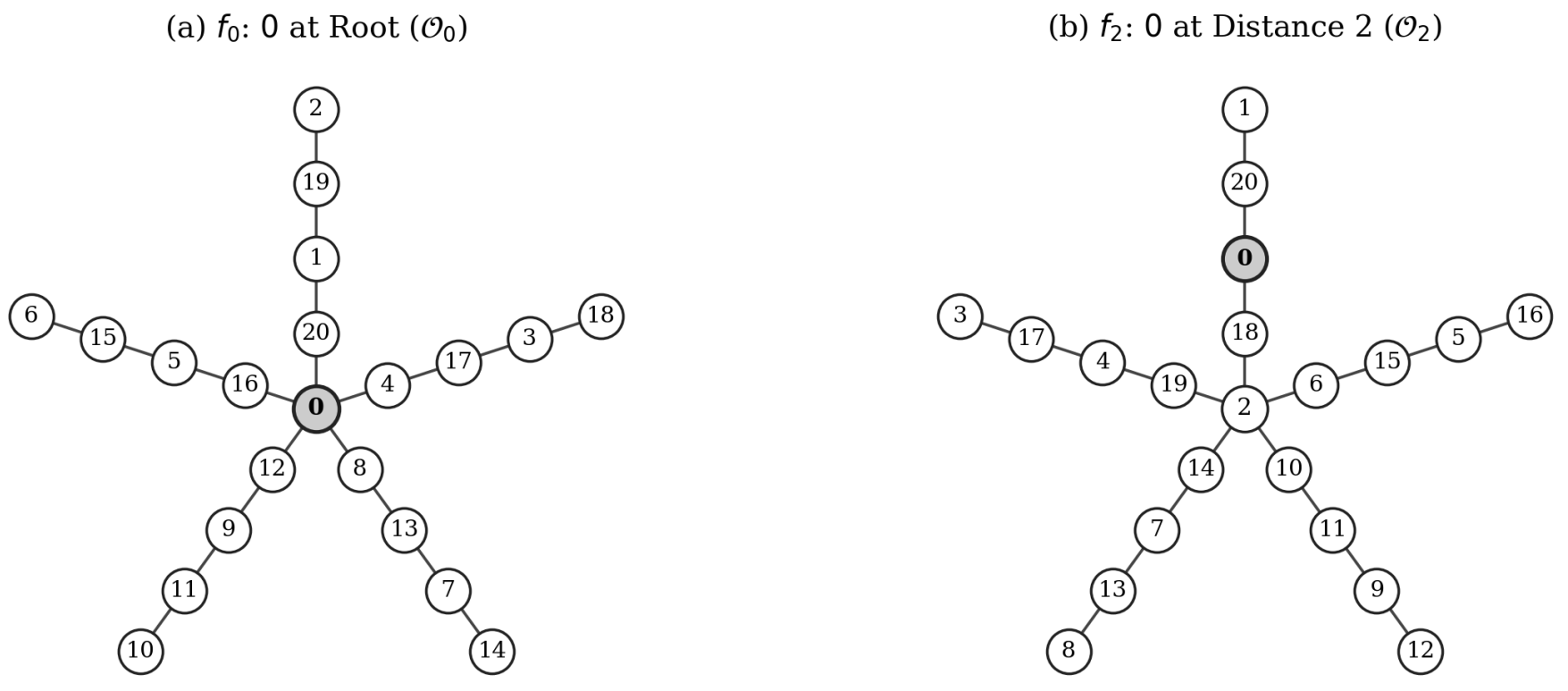

Figure 1: Graceful labelings of $T(5, 4)$ placing 0 (shaded) at the root $O_0$ and at distance 2 ($O_2$).

## References

1. F. R. K. Chung and F. K. Hwang, *Rotatability of trees*, Congr. Numer. **33** (1981), 239–250.

1. S. W. Golomb, *How to number a graph*, in *Graph Theory and Computing* (R. C. Read, ed.), Academic Press, New York, 1972, pp. 23–37.
2. R. I. Rofa, *A graceful algebraic function labelling of rooted symmetric trees*, arXiv:2109.09511 [math.CO], 2021.
3. R. I. Rofa, *0-rotatability of classes of rooted symmetric trees. Are rooted symmetric trees 0-rotatable?*, arXiv: 2312.16235 [math.CO], 2023.
4. A. Rosa, *On certain valuations of the vertices of a graph*, in *Theory of Graphs (Internat. Sympos., Rome, 1966)*, Gordon and Breach, New York, 1967, pp. 349–355.